\documentclass[12pt]{amsart} 
\usepackage{amsfonts,graphics,amsmath,amsthm,amsfonts,amscd,amssymb,amsmath,latexsym,multicol}
\usepackage{epsfig}
\usepackage{flafter}

\input diagrams

\makeatletter

\def\jobis#1{FF\fi
  \def\predicate{#1}%
  \edef\predicate{\expandafter\strip@prefix\meaning\predicate}%
  \edef\job{\jobname}%
  \ifx\job\predicate
}

\makeatother

\if\jobis{proposal}%
 \def\try{subsection}%
\else
  \def\try{section}%
\fi

\theoremstyle{plain}
\newtheorem{theorem}{Theorem}[\try]
\newtheorem{corollary}[theorem]{Corollary}
\newtheorem{lemma}[theorem]{Lemma}
\newtheorem{claim}[theorem]{Claim}

\newtheorem{definition-lemma}[theorem]{Definition-Lemma}

\newtheorem{definition}[theorem]{Definition}

\newtheorem{conjecture}[theorem]{Conjecture}

\def\ideal#1.{I_{#1}}
\def\ring#1.{\mathcal {O}_{#1}}
\def\fring#1.{\hat{\mathcal {O}}_{#1}}
\def\proj#1.{\mathbb {P}(#1)}
\def\pr #1.{\mathbb {P}^{#1}}
\def\dpr #1.{\hat{\mathbb {P}}^{#1}}
\def\af #1.{\mathbb A^{#1}}
\def\Hz #1.{\mathbb F_{#1}}
\def\Hbz #1.{\overline{\mathbb F}_{#1}}
\def\fb#1.{\underset #1 {\times}}
\def\rest#1.{\underset {\ \ring #1.} \to \otimes}
\def\au#1.{\operatorname {Aut}\,(#1)}
\def\deg#1.{\operatorname {deg } (#1)}
\def\pic#1.{\operatorname {Pic}\,(#1)}
\def\pico#1.{\operatorname{Pic}^0(#1)}
\def\picg#1.{\operatorname {Pic}^G(#1)}
\def\ner#1.{NS (#1)}
\def\rdown#1.{\llcorner#1\lrcorner}
\def\rfdown#1.{\lfloor{#1}\rfloor}
\def\rup#1.{\ulcorner{#1}\urcorner}
\def\rcup#1.{\lceil{#1}\rceil}
\def\cone#1.{\operatorname {NE}(#1)}
\def\ccone#1.{\overline{\operatorname {NE}}(#1)}
\def\coef#1.{\frac{(#1-1)}{#1}}
\def\vit#1.{D_{\langle #1 \rangle}}
\def\mm#1.{\overline {M}_{0,#1}}
\def\H1#1.{H^1(#1,{\ring #1.})}
\def\ac#1.{\overline {\mathbb F}_{#1}}

\def\adj#1.{\frac {#1-1}{#1}}
\def\spn#1.{\overline{#1}}
\def\pek#1.#2.{\Cal P^{#1}(#2)}
\def\plk#1.#2.{\Cal P^{\leq #1}(#2)}
\def\ev#1.{\operatorname{ev_{#1}}}
\def\ilist#1.{{#1}_1,{#1}_2,\dots}
\def\bminv#1.{(\nu_1,s_1;\nu_2,s_2;\dots ;\nu_{#1},s_{#1};\nu_{r+1})}
\def\zinv#1.{(\nu_1,s_1;\nu_2,s_2;\dots ;\nu_{#1},s_{#1};0)}
\def\iinv#1.{(\nu_1,s_1;\nu_2,s_2;\dots ;\nu_{#1},s_{#1};\infty)}

\def\llist#1.#2.{{#1}_1,{#1}_2,\dots,{#1}_{#2}}
\def\lomitlist#1.#2.{{#1}_1,{#1}_2,\dots,\hat {{#1}_i}, \dots, {#1}_{#2}}
\def\lomitlistz#1.#2.{{#1}_0,{#1}_1,\dots,\hat {{#1}_i}, \dots, {#1}_{#2}}
\def\loc#1.#2.{\Cal O_{#1,#2}}
\def\fderiv#1.#2.{\frac {\partial #1}{\partial #2}}
\def\deriv#1.#2.{\frac {d #1}{d #2}}
\def\map#1.#2.{#1 \longrightarrow #2}
\def\rmap#1.#2.{#1 \dasharrow #2}
\def\emb#1.#2.{#1 \hookrightarrow #2}
\def\non#1.#2.{\text {Spec }#1[\epsilon]/(\epsilon)^{#2}}
\def\Hi#1.#2.{\text {Hilb}^{#1}(#2)}
\def\sym#1.#2.{\operatorname {Sym}^{#1}(#2)}
\def\Hb#1.#2.{\text {Hilb}_{#1}(#2)}
\def\Hm#1.#2.{\Hom_{#1}(#2)}
\def\prd#1.#2.{{#1}_1\cdot {#1}_2\cdots {#1}_{#2}}
\def\Bl #1.#2.{\operatorname {Bl}_{#1}#2}
\def\pl #1.#2.{#1^{\otimes #2}}
\def\mgn#1.#2.{\overline {M}_{#1,#2}}
\def\ialist#1.#2.{{#1}_1 #2 {#1}_2, #2\dots}
\def\pair#1.#2.{\langle #1, #2\rangle}
\def\vandermonde#1.#2.{\left|
\begin{matrix}
1 & 1 & 1 & \dots & 1\\
{#1}_1 & {#1}_2 & {#1}_3 & \dots & {#1}_{#2}\\
{#1}_1^2 & {#1}_2^2 & {#1}_3^2 & \dots & {#1}_{#2}^2\\
\vdots & \vdots & \vdots & \ddots & \vdots\\
{#1}_1^{#2-1} & {#1}_2^{#2-1} & {#1}_2^{#2-1} & \dots & {#1}_{#2}^{#2-1}\\
\end{matrix}
\right|
}
\def\vandermondet#1.#2.{\left|
\begin{matrix}
1 & {#1}_1   & {#1}_1^2 & \dots & {#1}_1^{#2-1}\\
1 & {#1}_2   & {#1}_2^2 & \dots & {#1}_2^{#2-1}\\
1 & {#1}_3   & {#1}_3^2 & \dots & {#1}_3^{#2-1}\\
\vdots & \vdots & \vdots & \ddots & \vdots\\
1 & {#1}_{#2}& {#1}_{#2}^2 & \dots & {#1}_{#2}^{#2-1}\\
\end{matrix}
\right|
}
\def\gr#1.#2.{\mathbb{G}(#1,#2)}

\def\alist#1.#2.#3.{{#1}_1 #2 {#1}_2 #2\dots #2 {#1}_{#3}}
\def\zlist#1.#2.#3.{#1_0 #2 #1_1 #2\dots #2 #1_{#3}}
\def\lomitlist30#1.#2.#3.{{#1}_0,{#1}_1 #2 \dots #2\hat {{#1}_i} #2\dots #2 {#1}_{#3}}
\def\lmap#1.#2.#3.{#1 \overset{#2}{\longrightarrow} #3}
\def\mes#1.#2.#3.{#1 \longrightarrow #2 \longrightarrow #3}
\def\ses#1.#2.#3.{0\longrightarrow #1 \longrightarrow #2 \longrightarrow #3 \longrightarrow 0}
\def\les#1.#2.#3.{0\longrightarrow #1 \longrightarrow #2 \longrightarrow #3}
\def\res#1.#2.#3.{#1 \longrightarrow #2 \longrightarrow #3\longrightarrow 0}
\def\Hi#1.#2.#3.{\text {Hilb}^{#1}_{#2}(#3)}
\def\ten#1.#2.#3.{#1\underset {#2}{\otimes} #3}
\def\lomitlist30#1.#2.#3.{{#1}_0 #2 {#1}_1 #2 \dots #2 \hat {{#1}_i} #2 \dots #2 {#1}_{#3}}

\def\Hom{\operatorname{Hom}}

\def\Proj{\operatorname{Proj}}

\def\dim{\operatorname{dim}}

\def\deg{\operatorname{deg}}

\def\lcs{\operatorname{LCS}}

\def\mov{\operatorname{Mov}}

\def\rest{\operatorname{res}}

\def\e{\Cal E}

\def\e1{E_1}
\def\e2{E_2}

\def\ds{\displaystyle}

\def\bd{_{\bullet}}

\def\mapdown#1{\big\downarrow\rlap{$\vcenter
{\hbox{$\scriptstyle#1$}}$}}

\def\mapse#1{
{\vcenter{\hbox{$\mathop{\smash{\raise1pt\hbox{$\diagdown$}\!\lower7pt
\hbox{$\searrow$}}\vphantom{p}}\limits_{#1}\vphantom{\mapdown{}}$}}}}

\def\VR#1.{height#1pt&\omit&&\omit&&\omit&&\omit&&\omit&\cr}

\def\VRT#1.{height#1pt&\omit&&\omit&\cr}

\begin{document}
\title{On the existence of flips} 
\author{Christopher D. Hacon} 
\address{Department of Mathematics \\  
University of Utah\\  
155 South 1400 E\\
JWB 233\\
Salt Lake City, UT 84112, USA}
\email{hacon@math.utah.edu}
\author{James M\textsuperscript{c}Kernan} 
\address{Department of Mathematics\\ 
University of California at Santa Barbara\\ 
Santa Barbara, CA 93106, USA} 
\email{mckernan@math.ucsb.edu}

\thanks{The first author was partially supported by NSF research grant
  no: 0456363 and by a grant from the Sloan Foundation.  We would like
  to thank Adebisi Agboola, Florin Ambro, Paolo Cascini and Alessio
  Corti for many useful comments on the contents of this paper.}

\begin{abstract} Using the techniques of \cite{Siu98}, \cite{Kawamata99},
\cite{Shokurov93} and \cite{Shokurov03}, we prove that flips exist in dimension $n$, if
one assumes the termination of real flips in dimension $n-1$.
\end{abstract}

\maketitle

\section{Introduction}
\label{s_introduction}

The main result of this paper is:

\begin{theorem}\label{t_existence} Assume the real MMP in dimension $n-1$.  

 Then flips exist in dimension $n$.  
\end{theorem}

 Here are two consequences of this result:

\begin{corollary}\label{c_existence1} Assume termination of real flips in 
dimension $n-1$ and termination of flips in dimension $n$.  

 Then the MMP exists in dimension $n$.  
\end{corollary}

 As Shokurov has proved, \cite{Shokurov96}, the termination of real flips in dimension three, 
we get a new proof of the following result of Shokurov \cite{Shokurov03}:

\begin{corollary}\label{c_existence2} Flips exist in dimension four.  
\end{corollary}

Given a proper variety, it is natural to search for a good birational model.  An obvious,
albeit hard, first step is to pick a smooth projective model.  Unfortunately there are far
too many such models; indeed given any such, we can construct infinitely many more, simply
by virtue of successively blowing up smooth subvarieties.  To construct a unique model, or
at least cut down the examples to a manageable number, we have to impose some sort of
minimality on the birational model.

The choice of such a model depends on the global geometry of $X$.  One possibility is that
we can find a model on which the canonical divisor $K_X$ is nef, so that its intersection
with any curve is non-negative.  Conjecturally this is equivalent to the condition that
the Kodaira dimension of $X$ is non-negative, that is there are global pluricanonical forms.
Another possibility is that through any point of $X$ there passes a rational curve.  In
this case the canonical divisor is certainly negative on such a curve, and the best one
can hope for is that there is a fibration on which the anticanonical divisor is relatively
ample.  In other words we are searching for a birational model $X$ such that either
\begin{enumerate} 
\item $K_X$ is nef, in which case $X$ is called a \textit{minimal model}, or 
\item there is a fibration $\map X.Z.$ of relative Picard number one, such that $-K_X$ is
relatively ample; we call this a \textit{Mori fibre space}.
\end{enumerate} 

The minimal model program is an attempt to construct such a model step by step.  We start
with a smooth projective model $X$.  If $K_X$ is nef, then we have case (1).  Otherwise
the cone theorem guarantees the existence of a contraction morphism $f\colon \map X.Z.$ of
relative Picard number one, such that $-K_X$ is relatively ample.  If the dimension of $Z$
is less than the dimension of $X$, then we have case (2).  Otherwise $f$ is birational.
If $f$ is divisorial (that is the exceptional locus is a divisor), then we are free to
replace $X$ by $Z$ and continue this process.  Even though $Z$ may be singular, it is not
hard to prove that it is $\mathbb{Q}$-factorial, so that any Weil divisor is
$\mathbb{Q}$-Cartier (some multiple is Cartier), and that $X$ has terminal singularities.
In particular it still makes sense to ask whether $K_X$ is nef, and the cone theorem still
applies at this level of generality.  The tricky case is when $f$ is not divisorial, since
in this case it is not hard to show that no multiple of $K_Z$ is Cartier, and it no longer
even makes sense to ask if $K_Z$ is nef.  At this stage we have to construct the flip.

Let $f\colon\map X.Z.$ be a small projective morphism of normal varieties, so that $f$ is
birational but does not contract any divisors.  If $D$ is any integral Weil divisor such
that $-D$ is relatively ample the \textit{flip} of $D$, if it exists at all, is a
commutative diagram
$$
\begin{diagram}
X   &   & \rDashto   &  &  X'  \\
& \rdTo^f &  & \ldTo^{f'} &  \\
   &   & Z, &     &
\end{diagram}
$$
where $\rmap X.X'.$ is birational, and the strict transform $D'$ of $D$ is relatively
ample.  Note that $f'$ is unique, if it exists at all; indeed if we set
$$
\mathfrak{R}=R(X,D)=\bigoplus _{n\in\mathbb{N}} f_*\ring X.(nD),
$$
then 
$$
X'=\Proj _Z \mathfrak{R}.
$$
In particular the existence of the flip is equivalent to finite generation of the ring
$\mathfrak{R}$.

It is too much to expect the existence of general flips; we do however expect that flips
exist if $D=K_X+\Delta$ is kawamata log terminal.  Supposing that the flip of $D=K_X$
exists, we replace $X$ by $X'$ and continue.  Unfortunately this raises another issue, how
do we know that this process terminates?  It is clear that we cannot construct an
infinite sequence of divisorial contractions, since the Picard number drops by one after
every divisorial contraction, and the Picard number of $X$ is finite.  In other words, to
establish the existence of the MMP, it suffices to prove the existence and termination of
flips.  Thus \eqref{t_existence} reduces the existence of the MMP in dimension $n$, to
termination of flips in dimension $n$.

In the following two conjectures, we work with either the field $K=\mathbb{Q}$ or
$\mathbb{R}$.

\begin{conjecture}[Existence of Flips]\label{c_existence} Let $(X,\Delta)$ be a kawamata 
log terminal $\mathbb{Q}$-factorial pair of dimension $n$, where $\Delta$ is a
$K$-divisor.  Let $f\colon\map X.Z.$ be a flipping contraction, so that $-(K_X+\Delta)$
is relatively ample, and $f$ is a small contraction of relative Picard number one.

Then the flip of $f$ exists.  
\end{conjecture}

\begin{conjecture}[Termination of Flips]\label{c_termination} Let $(X,\Delta)$ be a kawamata 
log terminal $\mathbb{Q}$-factorial pair of dimension $n$, where $\Delta$ is a $K$-divisor.

 Then there is no infinite sequence of $(K_X+\Delta)$-flips.  
\end{conjecture}

 For us, the statement \lq\lq assuming the (real) MMP in dimension $n$\rq\rq\ means precisely
assuming \eqref{c_existence}$_{\mathbb{Q},n}$ and \eqref{c_termination}$_{\mathbb{R},n}$.  
In fact it is straightforward to see that \eqref{c_existence}$_{\mathbb{Q},n}$ implies 
\eqref{c_existence}$_{\mathbb{R},n}$, see for example the proof of \eqref{t_universal}.

Mori, in a landmark paper \cite{Mori88}, proved the existence of $3$-fold flips, when $X$
is terminal and $\Delta$ is empty.  Later on Shokurov \cite{Shokurov93} and Koll\'ar
\cite{Kollaretal} proved the existence of $3$-fold flips for kawamata log terminal pairs
$(X,\Delta)$, that is they proved \eqref{c_existence}$_3$.  Much more recently
\cite{Shokurov03}, Shokurov proved \eqref{c_existence}$_4$.  

Kawamata, \cite{Kawamata92c}, proved the termination of any sequence of threefold,
kawamata log terminal flips, that is he proved \eqref{c_termination}$_{\mathbb{Q},3}$.  As
previously pointed out, Shokurov proved, \cite{Shokurov96},
\eqref{c_termination}$_{\mathbb{R},3}$.  Further, Shokurov proved, \cite{Shokurov04}, that
\eqref{c_termination}$_{\mathbb{R},n}$, follows from two conjectures on the behaviour of
the log discrepancy of pairs $(X,\Delta)$ of dimension $n$ (namely acc for the set of log
discrepancies, whenever the coefficients of $\Delta$ are confined to belong to a set of
real numbers which satisfies dcc, and semicontinuity of the log discrepancy).  Finally,
Birkar, in a very recent preprint, \cite{Birkar05}, has reduced
\eqref{c_termination}$_{\mathbb{R},n}$, in the case when $K_X+\Delta$ has non-negative
Kodaira dimension, to acc for the log canonical threshold and the existence of the MMP in
dimension $n-1$.

 We also recall the abundance conjecture, 
\begin{conjecture}[Abundance]\label{c_abundance} Let $(X,\Delta)$ be a kawamata log terminal 
pair, where $K_X+\Delta$ is $\mathbb{Q}$-Cartier, and let $\pi\colon\map X.Z.$ be a proper
morphism, where $Z$ is affine and normal.

 If $K_X+\Delta$ is nef, then it is semiample.  
\end{conjecture}

Note that the three conjectures, existence and termination of flips, and abundance, are
the three most important conjectures in the MMP.  For example, Kawamata proved,
\cite{Kawamata85a}, that these three results imply additivity of Kodaira dimension.

 Our proof of \eqref{t_existence} follows the general strategy of \cite{Shokurov03}.  The
first key step was already established in \cite{Shokurov93}, see also \cite{Kollaretal}
and \cite{Fujino05}.  In fact it suffices to prove the existence of the flip for
$D=K_X+S+B$, see \eqref{t_reduction}, where $S$ has coefficient one, and $K_X+S+B$ is
purely log terminal.  The key point is that this allows us to restrict to $S$, and we can
try to apply induction.  By adjunction we may write
$$
(K_X+S+B)|_S=K_S+B',
$$
where $B'$ is effective and $K_S+B'$ is kawamata log terminal.  In fact, since we are
trying to prove finite generation of the ring $R$, the key point is to consider the
restriction maps
$$
\map {H^0(X,\ring X.(m(K_X+S+B)))}.{H^0(X,\ring S.(m(K_S+B')))}.,
$$
see \eqref{l_iff}.  Here and often elsewhere, we will assume that $Z$ is affine, so that
we can replace $f_*$ by $H^0$.  Now if these maps were surjective, we would be done by
induction.  Unfortunately this is too much to expect.  However we are able to prove, after
changing models, that something close to this does happen.

The starting point is to use the extension result proved in \cite{HM05b}, which in turn
builds on the work of Siu \cite{Siu98} and Kawamata \cite{Kawamata99}.  To apply this
result, we need to improve how $S$ sits inside $X$.  To this end, we pass to a resolution
$\map Y.X.$.  Let $T$ be the strict transform of $S$ and let $\Gamma$ be those divisors of
log discrepancy less than one.  Then by a generalisation of (3.17) of \cite{HM05b}, we can
extend sections from $T$ to $Y$, provided that we can find $G\in |m(K_Y+\Gamma)|$ which
does not contain any log canonical centre of $K_Y+\rup \Gamma.$.  

In fact if we blow up more, we can separate all of the components of $\Gamma$, except the
intersections with $T$.  In this case, the condition on $G$ becomes that it does not share
any components with $\Gamma$.  Thus, for each $m$ we are able to cancel common components,
and lift sections.  Putting all of this together, see \S \ref{s_extend} and \S \ref{s_limiting} for more details, we get a
sequence of divisors $\Theta\bd$ on $T$, such that
$$
i\Theta_i+j\Theta_j\leq (i+j)\Theta_{i+j},
$$
and it suffices to prove that this sequence stabilises, that is 
$$
\Theta_m=\Theta,
$$
is constant for $m$ sufficiently large and divisible.  To this end, we take the limit
$\Theta$ of this sequence.  Then $K_T+\Theta$ is kawamata log terminal, but in general
since $\Theta$ is a limit, it has real coefficients, rather than just rational.

Now there are two ways in which the sequence $\Theta\bd$ might vary.  By assumption each
$K_T+\Theta_m$ is big and so there is some model $\map T_m.T.$ on which the moving part of
$mk(K_T+\Theta_m)$ becomes semiample.  The problem is that the model $T_m$ depends on $m$.
This is obviously an issue of birational geometry, and can only occur in dimensions two
and higher.  To get around this, we need to run the real MMP, see \S \ref{s_competition}.
Thus replacing $T$ by a higher model, we may assume that the mobile part of some fixed
multple of $K_T+\Theta_m$ is in fact free, and the positive part of $K_X+\Theta$ is
semiample.

The second way in which which the sequence $\Theta\bd$ might vary is that we might get
freeness of the mobile part of $mk(K_T+\Theta_m)$ on the same model, but $\Theta_m$ is
still not constant.  There are plenty of such examples, even on the curve $\pr 1.$.
Fortunately Shokurov has already proved that this cannot happen, since the sequence
$\Theta\bd$ satisfies a subtle asymptoptic saturation property, see \S \ref{s_diophantine}
and \S \ref{s_saturation}.

Hopefully it is clear, from what we just said, the great debt our proof of
\eqref{t_existence} owes to the work of Kawamata, Siu and especially Shokurov.  The
material in \S \ref{s_extend} was inspired by the work of Siu \cite{Siu98} and Kawamata
\cite{Kawamata99} on deformation invariance of plurigenera, and lifting sections using
multiplier ideas sheaves.  On the other hand, a key step is to use the reduction to pl
flips, due to Shokurov contained in \cite{Shokurov93}.  Moreover, we use many of the
results and ideas contained in \cite{Shokurov03}, especially the notion of a saturated
algebra.  

Since the proof of \eqref{t_existence} is not very long, we have erred on the side of
making the proofs as complete as possible.  We also owe a great debt to the work of Ambro,
Fujino and especially Corti, who did such a good job of making the work of Shokurov more
accessible.  In particular much of the material contained in \S \ref{s_finite} and
\S \ref{s_competition}-\ref{s_saturation} is due to Shokurov, as well as some of the material
in the other sections, and we have followed the exposition of \cite{Corti05} and
\cite{Fujino05} quite closely.

\section{Notation and conventions}
\label{s_notation}

We work over the field of complex numbers $\mathbb{C}$.  A $\mathbb{Q}$-Cartier divisor
$D$ on a normal variety $X$ is \textit{nef} if $D\cdot C\geq 0$ for any curve $C\subset
X$.  We say that two $\mathbb{Q}$-divisors $D_1$, $D_2$ are $\mathbb{Q}$-linearly
equivalent ($D_1\sim _{\mathbb{Q}} D_2$) if there exists an integer $m>0$ such that $mD_i$
are linearly equivalent.  We say that a $\mathbb{Q}$-Weil divisor $D$ is big if we may
find an ample divisor $A$ and an effective divisor $B$, such that $D \sim _{\mathbb{Q}}
A+B$.  

A \textit{log pair} $(X,\Delta)$ is a normal variety $X$ and an effective
$\mathbb{Q}$-Weil divisor $\Delta$ such that $K_X+\Delta$ is $\mathbb{Q}$-Cartier.  We say
that a log pair $(X,\Delta)$ is \textit{log smooth}, if $X$ is smooth and the support of
$\Delta$ is a divisor with global normal crossings.  A projective morphism $g \colon\map
Y.X.$ is a \textit{log resolution} of the pair $(X,\Delta )$ if $Y$ is smooth and
$g^{-1}(\Delta )\cup \{\,\text {exceptional set of $g$}\,\}$ is a divisor with normal
crossings support.  We write $g^*(K_X +\Delta )=K_Y +\Gamma$ and $\Gamma =\sum a_i\Gamma
_i$ where $\Gamma _i$ are distinct reduced irreducible divisors.  The log discrepancy of
$\Gamma_i$ is $1-a_i$.  The \textit{locus of log canonical singularities of the pair
  $(X,\Delta)$}, denoted $\lcs(X,\Delta)$, is equal to the image of those components of
$\Gamma$ of coefficient at least one (equivalently log discrepancy at most zero).  The
pair $(X,\Delta )$ is \textit{kawamata log terminal} if for every (equivalently for one)
log resolution $g\colon\map Y.X.$ as above, the coefficients of $\Gamma$ are strictly less
than one, that is $a_i<1$ for all $i$.  Equivalently, the pair $(X,\Delta)$ is kawamata
log terminal if the locus of log canonical singularities is empty.  We say that the pair
$(X,\Delta)$ is \textit{purely log terminal} if the log discrepancy of any exceptional
divisor is greater than zero.

 We will also write
$$
K_Y+\Gamma=g^*(K_X+\Delta)+E,
$$
where $\Gamma$ and $E$ are effective, with no common components, and $E$ is $g$-exceptional.
Note that this decomposition is unique.   

Note that the group of Weil divisors with rational or real coefficients forms a vector
space, with a canonical basis given by the prime divisors.  If $A$ and $B$ are two $\mathbb{R}$-divisors,
then we let $(A,B]$ denote the line segment
$$
\{\, \lambda A+\mu B\,|\, \lambda+\mu=1,\lambda>0,\mu\geq 0 \,\}.
$$

Given an $\mathbb{R}$-divisor, $\|D\|$ denotes the sup norm with respect to this basis.
We say that $D'$ is \textit{sufficiently close to $D$} if there is a finite dimensional
vector space $V$ such that $D$ and $D'\in V$ and $D'$ belongs to a sufficiently small ball
of radius $\delta>0$ about $D$,
$$
\| D-D'\|<\delta.
$$

We recall some definitions involving divisors with real coefficients:

\begin{definition}\label{d_divisor} Let $X$ be a variety.  
\begin{enumerate} 
\item An \textbf{$\mathbb{R}$-Weil divisor} $D$ is an $\mathbb{R}$-linear combination of
prime divisors.
\item Two $\mathbb{R}$-divisors $D$ and $D'$ are \textbf{$\mathbb{R}$-linearly equivalent} if their
difference is an $\mathbb{R}$-linear combination of principal divisors.
\item An \textbf{$\mathbb{R}$-Cartier divisor} $D$ is an $\mathbb{R}$-linear combination of
Cartier divisors.
\item An $\mathbb{R}$-Cartier divisor $D$ is \textbf{ample} if it is strictly positive on the cone of 
curves minus the origin.  
\item An $\mathbb{R}$-divisor $D$ is \textbf{effective} if it is a positive real linear
combination of prime divisors.
\item An $\mathbb{R}$-Cartier divisor $D$ is \textbf{big} if it is the sum of an ample divisor 
and an effective divisor.   
\item An $\mathbb{R}$-Cartier divisor $D$ is \textbf{semiample} if there is a contraction 
$\pi\colon\map X.Y.$ such that $D$ is linearly equivalent to the pullback of an ample divisor.  
\end{enumerate} 
\end{definition}

Note that we may pullback $\mathbb{R}$-Cartier divisors, so that we may define the various
flavours of log terminal and log canonical in the obvious way.  

\section{Generalities on Finite generation}
\label{s_finite}

In this section we give some of the basic definitions and results concerning finite
generation; we only include the proofs for completeness.

We fix some notation.  Let $f\colon\map X.Z.$ be a projective morphism of
normal varieties, where $Z$ is affine.  Let $A$ be the coordinate ring of $Z$.  

\begin{definition-lemma}\label{d_truncation} Let $\mathfrak{R}$ be any graded $A$-algebra. 
A \textbf{truncation} of $\mathfrak{R}$ is any $A$-algebra of the form
$$
\mathfrak{R}_{(d)}=\bigoplus_{m\in\mathbb{N}}\mathfrak{R}_{md},
$$
for a positive integer $d$.  

 Then $\mathfrak{R}$ is finitely generated iff there is a positive integer $d$ such that $\mathfrak{R}_{(d)}$ 
is finitely generated.  
\end{definition-lemma}
\begin{proof} Suppose that $\mathfrak{R}$ is finitely generated.  The cyclic group $\mathbb{Z}_d$ 
acts in an obvious way on $\mathfrak{R}$, and under this action $\mathfrak{R}_{(d)}$ is the algebra of
invariants.  Thus $\mathfrak{R}_{(d)}$ is finitely generated by Noether's Theorem, which says that
the ring of invariants of a finitely generated ring, under the action of a finite group,
is finitely generated.  

 Now suppose that $\mathfrak{R}_{(d)}$ is finitely generated.  Let $f\in \mathfrak{R}$.  Then $f$ is a root of the 
monic polynomial
$$
x^d-f^d\in \mathfrak{R}_{(d)}[x].
$$
In particular $\mathfrak{R}$ is integral over $\mathfrak{R}_{(d)}$ and the result follows by Noether's
Theorem on the finiteness of the integral closure.
\end{proof}

 We are interested in finite generation of the following algebras:
\begin{definition}\label{d_divisorial} Let $B$ be an integral Weil divisor on $X$.  We
call any $\ring Z.$-algebra of the form 
$$
\bigoplus_{m\in\mathbb{N}}f_*\ring X.(mB),
$$
\textbf{divisorial}.
\end{definition}

 In particular we have:

\begin{lemma}\label{l_equivalence} Let $X$ be a normal variety and let 
$\mathfrak{R}$ and $\mathfrak{R}'$ be two divisorial algebras associated to divisors $D$ and 
$D'$.  

 If $aD\sim a'D'$ then $\mathfrak{R}$ is finitely generated iff $\mathfrak{R}'$ is finitely generated. 
\end{lemma}
\begin{proof} Clear, since $\mathfrak{R}$ and $\mathfrak{R}'$ have the same truncation.  
\end{proof}

  We want to restrict a divisorial algebra to a prime divisor $S$:

\begin{definition}\label{d_restricted} Let $\mathfrak{R}$ be the divisorial algebra 
associated to the divisor $B$.  The \textbf{restricted algebra} $\mathfrak{R}_S$ 
is the image 
$$
\map {\bigoplus_{m\in\mathbb{N}}f_*\ring X.(mB)}.{\bigoplus_{m\in\mathbb{N}}f_*\ring S.(mB)}.,
$$
under the obvious restriction map.  
\end{definition}

\begin{lemma}\label{l_iff} If the algebra $\mathfrak{R}$ is finitely generated then so is the 
restricted algebra.  Conversely, if $S$ is linearly equivalent to a multiple of $B$, where
$B$ is an effective divisor which does not contain $S$ and the restricted algebra is
finitely generated, then so is $\mathfrak{R}$.
\end{lemma}
\begin{proof} Since by definition there is a surjective homomorphism 
$$
\phi\colon\map \mathfrak{R}.\mathfrak{R}_S.,
$$
it follows that if $\mathfrak{R}$ is finitely generated then so is $\mathfrak{R}_S$.  

 Now suppose that $S\sim bB$.  Passing to a truncation, we may assume that $b=1$.  
We can identify the space of sections of 
$$
\mathfrak{R}_m=f_*\ring X.(mB)
$$
with rational functions $g$, such that 
$$
(g)+mB\geq 0.
$$

Let $g_1\in \mathfrak{R}_1$ be the rational function such that 
$$
(g_1)+B=S.
$$
Suppose we have $g\in \mathfrak{R}_m$, with $\phi(g)=0$.  Then the support of 
$$
(g)+mB,
$$
contains $S$, so that we may write
$$
(g)+mB=S+S',
$$
where $S'$ is effective.  But then 
$$
(g)+mB=(g_1)+B+S',
$$
so that 
$$
(g/g_1)+(m-1)B=S'.
$$
But this says exactly that $g/g_1\in \mathfrak{R}_{m-1}$, so that the kernel of $\phi$ is precisely
the principal ideal generated by $g_1$.  But then if $\mathfrak{R}_S$ is finitely generated, 
it is clear that $\mathfrak{R}$ is finitely generated.  
\end{proof}

\begin{definition}\label{d_additive} We say that a sequence of $\mathbb{R}$-divisors
$B\bd$ is \textbf{additive} if 
$$
B_i+B_j\leq B_{i+j},
$$
we say that it is \textbf{convex} if 
$$
\frac i{i+j}B_i+\frac j{i+j}B_j\leq B_{i+j},
$$
and we say that it is \textbf{bounded} if there is a divisor $B$ such that
$$
B_i\leq B.
$$
\end{definition}

Since the maps in \eqref{d_restricted} are not in general surjective, the restricted
algebra is not necessarily divisorial.  However we will be able to show that it is of the
following form:

\begin{definition}\label{d_geometric} Any $\ring Z.$-algebra of the form 
$$
\bigoplus_{m\in\mathbb{N}}f_*\ring X.(B_m),
$$
where $B\bd$ is an additive sequence of integral Weil divisors, will be called
\textbf{geometric}.
\end{definition}

We are interested in giving necessary and sufficient conditions for a divisorial or more
generally a geometric algebra to be finitely generated.

\begin{definition}\label{d_moving} Let $B$ be an integral divisor on $X$.  Let $F$ be the 
\textbf{fixed part} of the linear system $|B|$, and set $M=B-F$.  We may write
$$
|B|=|M|+F,
$$
We call $M=\mov B$ the \textbf{mobile part} of $B$ and we call $B=M+F$ the
\textbf{decomposition} of $B$ into its mobile and fixed part.  We say that a divisor is
\textbf{mobile} if the fixed part is empty.
\end{definition}

\begin{definition}\label{d_characteristic} Let $\mathfrak{R}$ be the geometric algebra associated to the convex
sequence $B\bd$.  Let 
$$
B_m=M_m+F_m,
$$
be the decomposition of $B_m$ into its mobile and fixed parts.  The sequence of divisors
$M\bd$ is called the \textbf{mobile sequence} and the sequence of $\mathbb{Q}$-divisors
$D\bd$ given by
$$
D_i=\frac {M_i}i,
$$
is called the \textbf{characteristic sequence}.  

 We say that $\mathfrak{R}$ is \textbf{free} if $M_m$ is base point free, for every $m$.  
\end{definition}

Clearly the mobile sequence is additive and the characteristic sequence is convex.  The
key point is that finite generation of a divisorial algebra only depends on the mobile
part in each degree, even up to a birational map:

\begin{lemma}\label{l_birational} Let $\mathfrak{R}$ be a geometric algebra asociated to the 
convex sequence $B\bd$.  Let $g\colon\map Y.X.$ be any birational morphism and let $\mathfrak{R}'$ be
the geometric algebra on $Y$ associated to a convex sequence $B'\bd$.

If the mobile part of $g^*B_i$ is equal to the mobile part of $B_i'$ then $\mathfrak{R}$ is finitely
generated iff $\mathfrak{R}'$ is finitely generated.  
\end{lemma}
\begin{proof} Clear. 
\end{proof}

\begin{lemma}\label{l_easy} Let $\mathfrak{R}$ be a free geometric algebra and let $D$ be the 
limit of the characteristic sequence.

 If $D=D_k$ for some positive integer $k$ then $\mathfrak{R}$ is finitely generated.  
\end{lemma}
\begin{proof} Passing to a truncation, we may assume that $D=D_1$.  But then
$$
mD=mD_1=mM_1\leq M_m=mD_m\leq mD, 
$$
and so $D=D_m$, for all positive integers $m$.  Let $h\colon\map X.W.$ be the contraction
over $Z$ associated to $M_1$, so that $M_1=h^*H$, for some very ample divisor on $W$.  We
have $g^*M_m=mg^*M_1=h^*(mH)$ and so the algebra $\mathfrak{R}$ is nothing more than the coordinate
ring of $W$ under the embedding of $W$ in $\pr n.$ given by $H$, which is easily seen to
be finitely generated by Serre vanishing.
\end{proof}

\section{Reduction to pl flips and finite generation}
\label{s_reduction}

We recall the definition of a pl flipping contraction:
\begin{definition}\label{d_pl-flip} We call a morphism $f\colon\map X.Z.$ or 
normal varieties, where $Z$ is affine, a \textbf{pl flipping contraction} if 
\begin{enumerate}
\item $f$ is a small birational contraction of relative Picard number one,
\item $X$ is $\mathbb{Q}$-factorial, 
\item $K_X+\Delta$ is purely log terminal, where $S=\rdown \Delta.$ is irreducible, and 
\item $-(K_X+\Delta)$ and $-S$ are ample. 
\end{enumerate}
\end{definition}

 Shokurov, \cite{Shokurov93}, see also \cite{Kollaretal} and \cite{Fujino05}, has shown:
\begin{theorem}[Shokurov]\label{t_reduction} To prove \eqref{t_existence} it suffices to
construct the flip of a pl flipping contraction.  
\end{theorem}

 The aim of the rest of the paper is to prove:

\begin{theorem}\label{t_restricted} Let $(X,\Delta)$ be a log pair of dimension $n$
and let $f\colon\map X.Z.$ be a morphism, where $Z$ is affine and normal.  Let $k$ be a
positive integer such that $D=k(K_X+\Delta)$ is Cartier, and let $\mathfrak{R}$ be the
divisorial algebra associated to $D$.  Assume that
\begin{enumerate} 
\item $K_X+\Delta$ is purely log terminal,   
\item $S=\rdown \Delta.$ is irreducible, 
\item there is a divisor $G\in |D|$, such that $S$ is not contained in the support of
$G$,
\item $\Delta-S\sim _{\mathbb{Q}}A+B$, where $A$ is ample and $B$ is an effective divisor, 
whose support does not contain $S$, and 
\item $-(K_X+\Delta)$ is ample.
\end{enumerate}

  If the real MMP holds in dimension $n-1$ then the restricted algebra $\mathfrak{R}_S$ is
finitely generated.
\end{theorem}

Note that the only interesting case of \eqref{t_restricted} is when $f$ is birational,
since otherwise the condition that $-(K_X+\Delta)$ is ample implies that
$\kappa(X,K_X+\Delta)=-\infty$.

We note that to prove \eqref{t_existence}, it is sufficient to prove \eqref{t_restricted}:

\begin{lemma}\label{l_implies} \eqref{t_restricted}$_n$ implies \eqref{t_existence}$_n$.  
\end{lemma}
\begin{proof} By \eqref{t_reduction} it suffices to prove the existence of pl flips.
Since $Z$ is affine and $f$ is small, it follows that $S$ is mobile.  By \eqref{l_iff} it
follows that it suffices to prove that the restricted algebra is finitely generated.
Hence it suffices to prove that a pl flip satisfies the hypothesis of
\eqref{t_restricted}.  Properties (1-2) and (5) are automatic and (3) follows as $S$ is
mobile.  $\Delta$ is automatically big, as $f$ is birational, and so $\Delta\sim
_{\mathbb{Q}} A+B$, where $A$ if ample, and $B$ is effective.  As $S$ is mobile, we may
assume that $B$ does not contain $S$.
\end{proof}

\section{Extending sections}
\label{s_extend}

The key idea of the proof of \eqref{t_existence} is to use the main result of \cite{HM05b}
to lift sections.  In this section, we show that we can improve this result, if we add
some hypotheses.  We recall some of the basic results about multiplier ideal sheaves.

\begin{definition}\label{d_multiplier} Let $(X,\Delta)$ be a log pair, where 
$X$ is smooth and let $\mu\colon\map Y.X.$ be a log resolution.  Suppose that 
we write
$$
K_Y+\Gamma=\mu^*(K_X+\Delta).
$$
The \textbf{multiplier ideal sheaf} of the log pair $(X,\Delta)$ is defined as
$$
\mathcal{J}(X,\Delta)=\mathcal{J}(\Delta)=\mu_*(-\rdown\Gamma.).
$$
\end{definition}

Note that the pair $(X,\Delta)$ is kawamata log terminal iff the multiplier ideal sheaf is
equal to $\ring X.$.  Another key property of a multiplier ideal sheaf is that it is 
independent of the log resolution.  Multiplier ideal sheaves have the following basic
property, see (2.2.1) of \cite{Takayama05}:

\begin{lemma}\label{l_basic} Let $(X,\Delta)$ be a kawamata log terminal pair, where 
$X$ is a smooth variety, and let $D$ be any divisor.  Let $f\colon\map X.Z.$ be any
projective morphism, where $Z$ is affine and normal.  Let $\sigma\in H^0(X,L)$ be any
section of a line bundle $L$, with zero locus $S\subset X$.

 If $D-S\leq \Delta$ then $\sigma\in H^0(X,(L\otimes \mathcal{J}(D)))$.  
\end{lemma}
\begin{proof} Let $g\colon\map Y.X.$ be a log resolution of the pair $(X,D+\Delta)$.  As
$S$ is integral
$$
\rdown g^*D.-g^*S\leq \rdown g^*\Delta.,
$$
and as the pair $(X,\Delta)$ is kawamata log terminal,
$$
K_{Y/X}-\rdown g^*\Delta.=-\rdown \Gamma.\geq 0.
$$
Thus
\begin{align*} 
g^*\sigma &\in H^0(Y,g^*L(-g^*S))\\ 
          &\subset H^0(Y,g^*L(-g^*S+K_{Y/X}-\rdown g^*\Delta.))\\ 
          &\subset H^0(Y,g^*L(K_{Y/X}-\rdown g^*D.)).\\ 
\end{align*} 
Pushing forward via $g$, we get 
$$
\sigma\in H^0(X,L\otimes \mathcal{J}(D)). 
$$
\end{proof}

We also have the following important vanishing result, which is an easy consequence of
Kawamata-Viehweg vanishing:

\begin{theorem}(Nadel Vanishing)\label{t_vanishing} Let $X$ be a smooth variety, let 
$\Delta$ be an effective divisor.  Let $f\colon\map X.Z.$ be any projective morphism and
let $N$ be any integral divisor such that $N-\Delta$ is relatively big and nef.

 Then
$$
R^if_*(\ring X.(K_X+N)\otimes \mathcal{J}(\Delta))=0, \qquad \text{for $i>0$.}
$$
\end{theorem}

 Here is the main result of this section:  

\begin{theorem}\label{t_lift} Let $(Y,\Gamma)$ be a smooth log pair and let $\pi\colon\map Y.Z.$ 
be a projective morphism, where $Z$ is normal and affine.  Let $m$ be a positive integer,
and let $L$ be any line bundle on $X$, such that $c_1(L)\sim _{\mathbb{Q}} m(K_Y+\Gamma)$.
Assume that
\begin{enumerate}
\item $(Y,\Gamma)$ is purely log terminal, 
\item $T=\rdown \Gamma.$ is irreducible, 
\item $\Gamma-T\sim_{\mathbb{Q}}A+B$, where $A$ is ample and $B$ is an effective divisor, 
which does not contain $T$.  
\end{enumerate}
Let $\Delta=(\Gamma-T)|_T$, so that 
$$
(K_Y+\Gamma)|_T=K_T+\Delta.
$$
Suppose that there is an effective divisor $H$, which does not contain $T$, such that for
every sufficiently divisible positive integer $l$, the natural homomorphism
$$
\map H^0(Y,L^{\otimes l}(H)).H^0(T,L^{\otimes l}(H)|_T).,
$$
contains the image of $H^0(T,L^{\otimes l}|_T)$, considered as a subspace of
$H^0(T,L^{\otimes l}(H)|_T)$ by the inclusion induced by $H$.  

Then the natural restriction homomorphism 
$$
\map H^0(Y,L).H^0(T,L|_T).,
$$
is surjective.  
\end{theorem} 
\begin{proof} As $K_Y+T+(1-\epsilon)\Gamma+\epsilon A+\epsilon B$ is purely log terminal
for any $\epsilon>0$ sufficiently small, replacing $A$ by $\epsilon A$ and $B$ by
$\epsilon B+(1-\epsilon)\Gamma$, we may assume that $K_Y+\Gamma=K_Y+T+A+B$ is purely log
terminal.  

We let primes denote restriction to $T$, so that, for example, $H'=H|_T$.  Fix a non-zero
section
$$
\sigma \in H^0(T,L').
$$
Let $S$ be the zero locus of $\sigma$.  By assumption, we may find a divisor $G_l\sim l c_1(L)+H$,
such that
$$
G_l'=lS+H',
$$
If we set
$$
N=c_1(L)-K_Y-T\qquad \text{and} \qquad \Theta =\frac{m-1}{ml}G_l+B,
$$
then 
$$
N\sim_{\mathbb{Q}}(m-1)(K_Y+\Gamma)+A+B.  
$$
Since 
$$
N-\Theta\sim _\mathbb{Q} A+B-\frac{m-1}{ml} H-B=A-\frac {m-1}{ml} H,
$$
is ample for $l$ sufficiently large, it follows that,
$$
H^1(Y,L(-T)\otimes \mathcal{J}(\Theta))=H^1(Y,\ring Y.(K_Y+N)\otimes \mathcal{J}(\Theta))=0,
$$
by Nadel vanishing \eqref{t_vanishing}, so that
$$
\map {H^0(Y,L\otimes \mathcal{J} (\Theta))}.{H^0(T,L'\otimes \mathcal{J}(\Theta))}.,
$$
is surjective.  Now 
\begin{align*} 
\Theta'-S &= B'+\frac{m-1}{ml}(lS+H')-S\\ 
            &\leq B'+\frac{m-1}{ml}H'.\\ 
\end{align*} 
Since $(Y,T+A+B)$ is purely log terminal, $(T,B')$ is kawamata log terminal,
and so 
$$
(T,B'+\frac{m-1}{ml}H'),
$$ 
is kawamata log terminal for $l$ sufficiently large.  But then 
$$
\sigma\in H^0(T,L'\otimes \mathcal{J}(\Theta'))\subset H^0(T,L'\otimes \mathcal{J}(\Theta)),
$$
by \eqref{l_basic}.  
\end{proof}

\section{Limiting algebras}
\label{s_limiting}

 To state the main result of this section, we need a:

\begin{definition}\label{d_limiting} We say that a geometric algebra $\mathfrak{R}$, given by an 
additive sequence $B\bd$, is \textbf{limiting}, if there are $\mathbb{Q}$-divisors
$\Delta_m$ and a positive integer $k$ such that 
\begin{enumerate} 
\item $B_m=mk(K_X+\Delta_m)$, 
\item the limit $\Delta$ of the convex sequence $\Delta\bd$ exists, and   
\item $K_X+\Delta$ is kawamata log terminal.
\end{enumerate} 
\end{definition}

\begin{theorem}\label{t_reduce} Let $(X,\Delta)$ be a log pair of dimension $n$
and let $f\colon\map X.Z.$ be a morphism, where $Z$ is affine and normal.  Let $k$ be any
positive integer such that $D=k(K_X+\Delta)$ is Cartier.  Assume that
\begin{enumerate} 
\item $K_X+\Delta$ is purely log terminal,   
\item $S=\rdown \Delta.$ is irreducible, 
\item there is a positive integer $m_0$ and a divisor $G_0\in |m_0D|$, such that 
$S$ is not contained in the support of $G_0$, and  
\item $\Delta-S\sim _{\mathbb{Q}}A+B$, where $A$ is ample and $B$ is an effective divisor,
whose support does not contain $S$.
\end{enumerate} 

 Then there is a log resolution $g\colon\map Y.X.$ with the following properties.  Suppose
that we write
$$
K_Y+\Gamma=g^*(K_X+\Delta)+E,
$$
where $\Gamma$ and $E$ are effective, with no common components and $E$ is exceptional.
Let $T$ be the strict transform of $S$ and let $\pi$ the composition of $f$ and $g$.  Let
$\mathfrak{R}$ be the divisorial algebra associated to $G=k(K_Y+\Gamma)$.

Then the restricted algebra $\mathfrak{R}_T$ is a limiting algebra, given by a convex
sequence $\Theta\bd$.
\end{theorem}

Note that the hypotheses of \eqref{t_reduce} are simply those of \eqref{t_restricted},
excluding (5) of \eqref{t_restricted}, and the hypothesis that the MMP holds.  To prove
\eqref{t_reduce}, we are going to apply \eqref{t_lift}.  The idea will be to start with
the main result (3.17) of \cite{HM05b}, which we state in a convenient form:

\begin{theorem}\label{t_extend} Let $(Y,\Gamma)$ be a smooth log pair, and let $\pi\colon\map Y.Z.$
be a projective morphism, where $Z$ is normal and affine.  Let $H=lA$, where $l$ is a
sufficiently large and divisible positive integer and $A$ is very ample.  Let $m$ be any
positive integer such that $m(K_Y+\Gamma)$ is Cartier.  Let $L$ be the line bundle $\ring
Y.(m(K_Y+\Gamma))$.  Assume that
\begin{enumerate}
\item $\Gamma$ contains $T$ with coefficient one, 
\item $(Y,\Gamma )$ is log canonical, and
\item there is a positive integer $m_0$ and a divisor $G_0\in |m_0(K_Y+\Gamma)|$ which 
does not contain any log canonical centre of $K_Y+\rup \Gamma.$.  
\end{enumerate}
Let $\Theta=(\Gamma-T)|_T$, so that 
$$
(K_Y+\Gamma)|_T=K_T+\Theta.
$$
Suppose that $H$ does not contain $T$.  Then the image of the natural homomorphism
$$
\map H^0(Y,L(H)).H^0(T,L(H)|_T).,
$$
contains the image of $H^0(T,L|_T)$, where $H^0(T,L|_T)$ is considered as a subspace of
$H^0(T,L(H)|_T)$ by the inclusion induced by $H|_T$.
\end{theorem}

Now to apply \eqref{t_extend}, the main point will be to change models and alter $\Gamma$,
so that property (3) holds.  To this end, we will need some results concerning
manipulation of log pairs.  Given a divisor $\Delta=\sum _i a_i\Delta_i$, we set
$$
\langle \Delta\rangle=\sum _i b_i \Delta_i \qquad \text{where} \qquad b_i = \begin{cases} a_i & \text {if $0<a_i<1$} \\
                                                                                           0  & \text{otherwise.} \\
\end{cases}
$$

\begin{lemma}\label{l_extend} Let $(Y,\Gamma)$ be a smooth log pair, and let $\pi\colon\map Y.Z.$
be a projective morphism, where $Z$ is normal and affine.  Let $m$ be any positive integer
such that $m(K_Y+\Gamma)$ is integral.  Let $L$ be the line bundle $\ring
Y.(m(K_Y+\Gamma))$.  Assume that
\begin{enumerate}
\item $(Y,\Gamma)$ is purely log terminal, 
\item $T=\rdown\Gamma.$ is irreducible, 
\item no two components of $\langle\Gamma\rangle$ intersect, and
\item there is a divisor $G\in |m(K_Y+\Gamma)|$ such that $G$ and $\Gamma$ have no common
components.  
\end{enumerate}
Let $\Theta'=(\Gamma-T)|_T$, so that 
$$
(K_Y+\Gamma)|_T=K_T+\Theta'.
$$
Then we may find a $\mathbb{Q}$-divisor $0\leq \Theta\leq\Theta'$ such that the image of
the natural homomorphism
$$
\map H^0(Y,L).H^0(T,L|_T).,
$$
may be identified with $H^0(T,\ring T.(m(K_T+\Theta)))$, considered as a subspace of
$H^0(T,L|_T)$ by the inclusion induced by $m(\Theta'-\Theta)$.  
\end{lemma}
\begin{proof} Since no two components of $\langle\Gamma\rangle$ intersect, and $T$ is the
only component of coefficient one, the only possible log canonical centres of $K_Y+\rup
\Gamma.$ contained in $G$, are the components of $T\cap \langle\Gamma\rangle$.

It follows that there is a resolution $h\colon\map Y'.Y.$ of the base locus of
$m(K_Y+\Gamma)$, which is a sequence of smooth blow ups with centres equal to the
irreducible components of $T\cap \langle\Gamma\rangle$, with the following property.  

 We may write 
$$
K_{Y'}+\Gamma'=h^*(K_Y+\Gamma)+E,
$$ 
where $\Gamma'$ and $E$ are effective, with no common components, and $E$ is exceptional.
Note that $m(K_{Y'}+\Gamma')$ and $mE$ are integral and that $G'=h^*G+mE\in
|m(K_{Y'}+\Gamma')|$.  Let
$$
m(K_{Y'}+\Gamma')=N_m+G_m,
$$
be the decomposition of $m(K_{Y'}+\Gamma')$ into its moving and fixed parts.  Then the
base locus of $N_m$ does not contain any log canonical centre of $K_{Y'}+\rup \Gamma'.$.

Cancelling common components of $G_m$ and $\Gamma'$, we may therefore find divisors $G'_m$
and $\Gamma'_m$, with no common components, such that
$$
m(K_{Y'}+\Gamma'_m)=N_m+G'_m,
$$
is the decomposition of $m(K_{Y'}+\Gamma'_m)$ into its moving and fixed parts.  Let $T'$
be the strict transform of $T$.  Since $h$ is a composition of blow ups, with smooth
centres, which intersect $T$ in a divisor, in fact $h|_{T'}\colon\map T'.T.$ is an
isomorphism.  Set $L'=\ring Y'.(m( K_{Y'}+\Gamma'_m))$. 

Possibly replacing $kA$ by a linearly equivalent divisor, we may assume that $g^*A$ and
the strict transform of $A$ are equal.  Since $A$ is ample, there is an effective and
exceptional divisor $F$ such that $g^*A-F$ is ample.  In this case
$$
\Gamma'_m-T'\sim _{\mathbb{Q}} (g^*A-F)+(\Gamma'_m-T'-g^*A+F)=A'+B'.
$$  
As there is a natural identification
$$
H^0(Y,L)=H^0(Y',L'),
$$
we are thus free to replace the pair $(Y,L)$ by $(Y',L')$, so that, letting $\Theta
=(\Gamma '_m -T')|_{T'}$, the result follows by \eqref{t_extend} and \eqref{t_lift}.
\end{proof}

\begin{lemma}\label{l_disjoint} Let $(X,\Delta)$ be a log pair.  We may find a birational
projective morphism
$$
g\colon\map Y.X.,
$$ 
with the following properties.  Suppose that we write
$$
K_Y+\Gamma=g^*(K_X+\Delta)+E,
$$
where $\Gamma$ and $E$ are effective, with no common components, and $E$ is exceptional.  

Then no two components of $\langle \Gamma\rangle$ intersect.
\end{lemma}
\begin{proof} Passing to a log resolution, we may assume that the pair $(X,\Delta)$ has
global normal crossings.  We will construct $g$ as a sequence of blow ups of irreducible
components of the intersection of a collection of components of $\langle \Delta\rangle$.
Now if a collection of components of $\langle \Delta\rangle$ intersect, then certainly no
irreducible component of their intersection is contained in $\rdown\Delta.$.  Since
$(X,\Delta)$ has global normal crossings, it follows that we may as well replace $\Delta$
by $\langle \Delta\rangle$.  Thus we may assume that the coefficients of the components of
$\Delta$ are all less than one, so that the pair $(X,\Delta)$ is kawamata log terminal,
and our aim is to find $g$, so that no two components of $\Gamma$ intersect.

We proceed by induction on the maximum number $k$ of components of $\Delta$ which
intersect.  Since $(X,\Delta)$ has normal crossings, $k\leq n=\dim X$, and it suffices to
decrease $k$.  We now proceed by induction on the maximum sum $s$ of the coefficients of
$k$ components which intersect.  If we pick $r$ such that $r\Delta$ is integral, then $s$
is at least $k/r$ and $rs$ is an integer, so it suffices to decrease $s$.  We further
proceed by induction on the number $l$ of subvarieties $V$ which are the components of the
intersection of $k$ components of $\Delta$ whose coefficients sum to $s$.  We aim to
decrease $l$ by blowing up.

Suppose that we blow up $g\colon\map Y.X.$ along the intersection $V$ of $k$ components
$\llist \Delta.k.$ of $\Delta$, with coefficients $\llist a.k.$.  A simple calculation,
see for example (2.29) of \cite{KM98}, gives that the discrepancy of the exceptional
divisor $E$ is $(k-1)-s$, so that
$$
K_Y+\Gamma=g^*(K_X+\Delta)+(k-1-s)E. 
$$
If $k-1-s\geq 0$, then $E$ is not a component of $\Gamma$.  Otherwise $E$ is a
component of $\Gamma$, with coefficient $s+1-k$.  Let $\llist \Gamma.k.$ be the components
of $\Gamma'$ which are the strict transforms of $\llist \Delta.k.$.  Then there are $k$
subvarieties of $Y$ which dominate $V$, which are the intersection of $k$ components of
$\Gamma$, namely the intersection with $E$ of all but one of $\llist\Gamma.k.$.  If $a_k$
is the smallest coefficient, then the maximum sum of the coefficients of these
intersections is
$$
(s-a_k)+(s+1-k)=s+[(s-a_k)-(k-1)]<s,
$$
and so we have decreased $l$ by one.  \end{proof}

\begin{proof}[Proof of \eqref{t_reduce}] Let $g\colon\map Y.X.$ be any morphism, 
whose existence is guaranteed by \eqref{l_disjoint}.  We may write 
$$
K_Y+\Gamma=g^*(K_X+\Delta)+E,
$$
where $\Gamma$ and $E$ are effective, with no common components and $E$ is exceptional.
Since $k(K_X+\Delta)$ is Cartier, $k(K_Y+\Gamma)$ and $kE$ are integral.  Let $T$ be the
strict transform of $S$ and let $\pi$ the composition of $f$ and $g$.  Let
$$
mk(K_Y+\Gamma)=N_m+G_m,
$$
be the decomposition of $mk(K_Y+\Gamma)$ into its moving and fixed parts.  By assumption,
$T$ is not a component of $G_m$.  Possibly replacing $k$ by a multiple, we may assume that
$kA$ is very ample.  Possibly replacing $kA$ by a linearly equivalent divisor, we may
assume that $g^*A$ and the strict transform of $A$ are equal.

Cancelling common components of $\Gamma$ and $G_m$, we may find divisors
$T+g^*A\leq\Gamma_m\leq\Gamma$ and $G'_m$, with no common components, such that
$$
mk(K_Y+\Gamma_m)=N_m+G'_m.
$$
Set $\Theta'_m=(\Gamma_m-T)|_T$ and $\Theta=(\Gamma-T)|_T$.  Let $L=\ring
Y.(mk(K_Y+\Gamma_m))$.   By\eqref{l_extend}, there is a divisor $\Theta_m\leq\Theta'_m$,
such that the image of the natural homomorphism
$$
\map H^0(Y,L).H^0(T,L|_T).,
$$
is equal to $H^0(T,\ring T.(m(K_T+\Theta_m)))$, considered as a subspace of $H^0(T,L|_T)$
by the inclusion induced by $m(\Theta'_m-\Theta_m)$.  On the other hand, as
$\Theta_m\leq\Theta$, the limit $\Theta'$ of the sequence $\Theta\bd$ exists and
$K_T+\Theta'$ is kawamata log terminal.  
\end{proof}

\section{Real versus rational} 
\label{s_competition}

 Most of the ideas and a significant part of the proofs of the results in this section 
are contained in \cite{Shokurov96}.  We have only restated these results at the level 
of generality we need to prove \eqref{t_existence}.  

We will need a generalisation of the base point free theorem to the case of real divisors:

\begin{theorem}[Base Point Free Theorem]\label{t_base} Let $(X,\Delta)$ 
be a $\mathbb{Q}$-factorial kawamata log terminal pair, where $\Delta$ is a
$\mathbb{R}$-divisor.  Let $f\colon\map X.Z.$ be a projective morphism, where $Z$ is
affine and normal, and let $D$ be a nef $\mathbb{R}$-divisor, such that $aD-(K_X+\Delta)$
is nef and big, for some positive real number $a$.

 Then $D$ is semiample.  
\end{theorem}
\begin{proof} Replacing $D$ by $aD$ we may assume that $a=1$.  By assumption we may 
write 
$$
D-(K_X+\Delta)=A+E,
$$
where $A$ is ample and $E$ is effective.  Thus
$$
D-(K_X+\Delta+\epsilon E),
$$
is ample for all $\epsilon>0$.  Since the pair $(X,\Delta+\epsilon E)$ is kawamata log
terminal for $\epsilon$ small enough, replacing $\Delta$ by $\Delta+\epsilon E$, we may
assume that
$$
D-(K_X+\Delta),
$$
is ample.  Perturbing $\Delta$, we may therefore assume that $K_X+\Delta$ is
$\mathbb{Q}$-Cartier.

Let $F$ be the set of all elements $\alpha$ of the closed cone of curves on which $D$ is
zero.  Then $K_X+\Delta$ is negative on $F$.  Let $H$ be any ample divisor.  For every ray
$R=\mathbb{R}^+\alpha$ contained in $F$, there is an $\epsilon>0$ such that
$K_X+\Delta+\epsilon H$ is negative on $R$.  By compactness of a slice, it follows that
there is an $\epsilon>0$, such that $K_X+\Delta+\epsilon H$ is negative on the whole of
$F$.  It follows by the cone theorem that $F$ is the span of finitely many extremal rays
$\llist R.k.$, where each extremal ray $R_i$ is spanned by an integral curve $C_i$.  Let
$\llist D.k.$ be the prime components of $D$.  Consider the convex subset $\mathcal{P}$ of
$$
\{\, B=\sum _i d_i D_i \,|\, d_i\in\mathbb{R}\,\},
$$
consisting of all divisors $B$ such that $B$ is zero on $F$.  Then $\mathcal{P}$ is a
closed rational polyhedral cone.

In particular $D\in\mathcal{P}$ is a convex linear combination of divisors $B_i\in
\mathcal{P}\cap N_{\mathbb{Q}}\cap U$, where $U$ is any neighbourhood of $D$.  But if $U$ is
sufficiently small, then $B_i-(K_X+\Delta)$ is also ample.  Now if
$B_i=(B_i-(K_X+\Delta))+(K_X+\Delta)$ is not nef, then it must be negative on a
$(K_X+\Delta)$-extremal ray.  As the extremal rays of $K_X+\Delta$ are discrete in a
neighbourhood of $F$, it follows that $B_i$ is also nef if $U$ is sufficiently small.  By
the base point free theorem, it follows that each $B_i$ is semiample, so that $D$ is
semiample.  \end{proof}

\begin{theorem}\label{t_universal} Assume the real MMP in dimension $n$.  Let $(X,\Delta)$ 
be a kawamata log terminal pair of dimension $n$, such that $K_X+\Delta$ is
$\mathbb{R}$-Cartier and big.  Assume that there is a divisor $\Psi$ such that $K_X+\Psi$
is $\mathbb{Q}$-Cartier and kawamata log terminal.  Let $f\colon\map X.Z.$ be any proper
morphism, where $Z$ is normal and affine.  Fix a finite dimensional vector subspace $V$ of 
the space of $\mathbb{R}$-divisors containing $\Delta$. 
 
Then there are finitely many birational maps $\psi_i\colon\rmap X.W_i.$, $1\leq i\leq l$
over $Z$, such that for every divisor $\Theta\in V$ sufficiently close to $\Delta$, there
is an $1\leq i\leq l$ with the following properties:
\begin{enumerate} 
\item $\psi_i$ is the composition of a sequence of $(K_X+\Theta)$-negative divisorial
contractions and birational maps, which are isomorphisms in codimension two, 
\item $W_i$ is $\mathbb{Q}$-factorial, and
\item $K_{W_i}+\psi_{i*}\Theta$ is semiample.  
\end{enumerate} 

 Further there is a positive integer $k$ such that 
\begin{enumerate} 
\setcounter{enumi}{3}
\item if $r(K_X+\Theta)$ is integral then $kr(K_{W_i}+\psi_{i*}\Theta)$ is base point free.  
\end{enumerate} 
\end{theorem}
\begin{proof} As the property of being big is an open condition, we may
asssume that for any $\Theta\in V$ sufficiently close to $\Delta$, $K_X+\Theta$ is big. 

Suppose that we have established (1-3).  As $W_i$ is $\mathbb{Q}$-factorial, it follows
that the group of Weil divisors modulo Cartier divisors is a finite group.  Thus there is
a fixed positive integer $s_i$ such that if $r(K_X+\Theta)$ is integral, then
$s_ir(K_{W_i}+\psi_{i*}\Theta)$ is Cartier.  By Koll\'ar's effective base point free
theorem, \cite{Kollar93b}, there is then a positive integer $M$ such that
$Ms_ir(K_{W_i}+\psi_{i*}\Theta)$ is base point free.  If we set $k$ to be $Ms$, where $s$
is the maximum of the $s_i$, then this is (4).  Thus it suffices to prove (1-3).

Now if $\Theta$ is sufficiently close to $\Delta$, then $K_X+\Theta$ is big, so that by
\eqref{t_base} we may replace (3) by the weaker condition,
\begin{itemize} 
\item[(3$'$)] $K_{W_i}+\psi_{i*}\Theta$ is nef.  
\end{itemize} 

Thus it suffices to establish (1), (2) and (3$'$).  
 
Since we are assuming existence and termination of flips for $\mathbb{Q}$-divisors, we may
construct a log terminal model of $(X,\Psi)$.  As $(X,\Psi)$ is kawamata log terminal, the
log terminal model is small over $X$.  Thus passing to a log terminal model of $(X,\Psi)$,
we may assume that $X$ is $\mathbb{Q}$-factorial and that $f$ is projective.

Suppose that $K_X+\Delta$ is not nef.  Let $R$ be an extremal ray for $K_X+\Delta$.  $R$
is necessarily $(K_X+\Theta)$-negative, for any $\mathbb{Q}$-divisor $\Theta$ close enough
to $\Delta$.  By the cone and contraction theorems applied to $K_X+\Theta$, we can
contract $R$, $\psi\colon\map X.X'.$.  $\psi$ must be birational, as $K_X+\Delta$ is big.
If $\psi$ is divisorial (that is the exceptional locus is a divisor) then we replace the
pair $(X,\Delta)$ by the pair $(X',\psi_*\Delta)$.  If $\psi$ is small, then using
\eqref{c_existence}$_{\mathbb{Q},n}$, we know the flip of $K_X+\Theta$ exists.  But then
this is also the flip of $K_X+\Delta$, and so we can replace the pair $(X,\Delta)$ by the
flip.  Since we are assuming \eqref{c_termination}$_{\mathbb{R},n}$, and we can only make
finitely many divisorial contractions, we must eventually arrive at the case when
$K_X+\Delta$ is nef.

By \eqref{t_base} it follows that $K_X+\Delta$ is relatively semiample.  Let
$\psi\colon\map X.W.$ be the corresponding contraction over $Z$.  Then there is an ample
$\mathbb{R}$-divisor $H$ on $W$ such that $K_X+\Delta=\psi^*H$.  Thus if $\Theta$ is
sufficiently close to $\Delta$ and $K_X+\Theta$ is relatively nef over $W$, then
$$
K_X+\Theta=K_X+\Delta+(\Theta-\Delta)=\psi^*H+(\Theta-\Delta),
$$
is nef.  

Note that we may replace $Z$ by $W$, and use the fact that a divisor is relatively
generated iff it is locally base point free.  Thus replacing $Z$ by an open affine subset
of $W$, we may assume that $f$ is birational and $K_X+\Delta$ is $\mathbb{R}$-linearly
equivalent to zero.  Let $B$ be the closure in $V$ of a ball with radius $\delta$ centred
at $\Delta$.  If $\delta$ is sufficiently small, then for every $\Theta\in B$,
$K_X+\Theta$ is kawamata log terminal.  Pick $\Theta$ a point of the boundary of $B$.
Since $K_X+\Delta$ is $\mathbb{R}$-linearly equivalent to zero, note that for every curve
$C$,
$$
(K_X+\Theta)\cdot C<0 \qquad \text{iff} \qquad (K_X+\Theta')\cdot C<0, \forall \Theta'\in (\Delta,\Theta].
$$

In particular every step of the $(K_X+\Theta)$-MMP is a step of $(K_X+\Theta')$-MMP, for
every $\Theta'\in (\Delta,\Theta]$.  Since we are assuming existence and termination of
flips, we have a birational map $\psi\colon\rmap X.W.$ over $Z$, such that
$K_W+\psi_*\Theta$ is nef, and it is clear that $K_W+\psi_*\Theta'$ is nef, for every
$\Theta'\in (\Delta,\Theta]$.

At this point we want to proceed by induction on the dimension of $B$.  To this end, note
that as $B$ is compact and $\Delta$ is arbitrary, our result is equivalent to proving that
(3$'$) holds in $B$.  By what we just said, this is equivalent to proving that (3$'$)
holds on the boundary of $B$, which is a compact polyhedral cone (since we are working in
the sup norm) and we are done by induction on the dimension of $B$.
\end{proof}

 The key consequence of \eqref{t_universal} is:

\begin{corollary}\label{c_universal} Assume the real MMP in dimension $n$.  Let 
$(X,\Delta)$ be a kawamata log terminal $\mathbb{Q}$-factorial pair of dimension $n$,
where $K_X+\Delta$ is an $\mathbb{R}$-divisor.  Let $f\colon\map X.Z.$ be a contraction
morphism.   Let $r$ be a positive integer.  

 If $K_X+\Delta$ is relatively big, then there is a birational model $g\colon\map Y.X.$ and
a positive integer $k$, such that if $\pi\colon\map Y.Z.$ is the composition of $f$ and
$g$, then for every divisor $\Theta$ sufficiently close to $\Delta$, 
\begin{enumerate} 
\item if $r(K_X+\Theta)$ is integral, then the moving part of $g^*(rk(K_X+\Theta))$ is
base point free. 
\item If $\Theta\bd$ is a convex sequence of divisors with limit $\Theta$, such that
$mr(K_X+\Theta_m)$ is integral then the limit $D$ of the characteristic sequence $D\bd$
associated to $B_m=g^*(mrk(K_X+\Theta_m))$ is semiample.
\end{enumerate} 
\end{corollary}
\begin{proof} Let $\psi_i\colon\rmap X.W_i.$ be the models, whose existence is guaranteed by
\eqref{t_universal}, and let $g\colon\map Y.X.$ be any birational morphism which resolves
the indeterminancy of $\psi_i$, $1\leq i\leq l$.  Let $\phi_i\colon\map Y.W_i.$ be the
induced birational morphisms, so that we have commutative diagrams
$$
\begin{diagram}
 &  &   Y   &  & \\
 & \ldTo^g &  & \rdTo^{\phi_i} \\
 X &  &   \rDashto^{\psi_i}    &  & W_i.
\end{diagram}
$$

Let $\Theta$ be sufficiently close to $\Delta$.  Then for some $i$,
$K_{W_i}+\psi_{i*}\Theta$ is semiample.  Suppressing the index $i$, we may write
$$
g^*(K_X+\Theta)=\phi^*(K_W+\psi_*\Theta)+E+F,
$$
where for $E$ we sum over the common exceptional divisors of $g$ and $\phi$, and for $F$
we sum over the exceptional divisors of $\phi$ which are not exceptional for $g$ (by
assumption there are no exceptional divisors of $g$ which are not also exceptional for
$\phi$).  (1) of \eqref{t_universal} implies that $F$ is effective.  But then by
negativity of contraction, see (2.19) of \cite{Kollaretal}, $E$ is also effective.

Suppose that $r(K_X+\Theta)$ is integral.  Then the moving part of $g^*(rk(K_X+\Theta))$
is equal to the moving part of $\phi^*(rk(K_W+\psi_*\Theta))$, and we can apply
\eqref{t_universal} to conclude that there is a fixed $k$ such that the moving part of
$g^*(rk(K_X+\Theta))$ is base point free.  This is (1).

Now suppose that $\Theta\bd$ is a convex sequence with limit $\Theta$, such that
$mr(K_X+\Theta_m)$ is integral.  Let $M_m$ be the mobile part of
$g^*(mrk(K_X+\Theta_m))$.  Then, by what we have already said, $M_m$ is also the mobile
part of $\phi_i^*(mrk(K_W+\psi_{i*}\Theta_m))$.  Possibly passing to a subsequence, we may
assume that $i$ is constant, and in this case we suppress it.  It is then clear that the
limit $D$ of 
$$
D_m=\frac{M_m}m,
$$
is 
$$
\phi^*(rk(K_W+\psi_*\Theta)),
$$
so that $D$ is nef.  It follows that $D$ is semiample by \eqref{t_universal} (or indeed
\eqref{t_base}).  \end{proof}

\section{Diophantine Approximation}
\label{s_diophantine}

All of the results in this section are implicit in the work of Shokurov \cite{Shokurov03},
and we claim no originality.  In fact we have only taken Corti's excellent introduction to
Shokurov's work on the existence of flips and restated those results without the use of
b-divisors.

\begin{lemma}[Diophantine Approximation]\label{l_diophantine} Let $Y$ be a smooth
variety and let $\pi\colon\map Y.Z.$ be a projective morphism, where $Z$ is affine and
normal.  Let $D$ be a semiample divisor on $Y$.  Let $\epsilon>0$ be a positive rational
number.

 Then there is an integral divisor $M$ and a positive integer $m$ such that 
\begin{enumerate} 
\item $M$ is base point free, 
\item $\ds \|mD-M\| < \epsilon$, and 
\item If $mD\geq M$ then $mD=M$.  
\end{enumerate} 
\end{lemma}
\begin{proof} If $D$ is rational, then pick $m$ such that $mD$ is integral and set 
$M=mD$.  Thus we may suppose that $D$ is not rational.  

Let $N_{\mathbb{Z}}$ be the lattice spanned by the components $G_j$ of $D$, and let
$N_{\mathbb{Q}}$ and $N_{\mathbb{R}}$ be the corresponding vector spaces.  Since $D$ is
semiample and $\pi$ is projective, we may pick a basis $\{P_k\}$ of $N_{\mathbb{Q}}$,
where each $P_k$ is base point free, and $D$ belongs to the cone
$$
\mathcal{P}=\sum\mathbb{R}_+[P_k]\subset \mathbb{R}_+[G_k]=\mathcal{G}.  
$$
Let $v\in N_{\mathbb{R}}$ be the vector corresponding to $D$.  Let $A$ be the cyclic 
subgroup of the torus 
$$
\frac{N_{\mathbb{R}}}{N_{\mathbb{Z}}},
$$
generated by the image of $v$.  Let $\bar A$ be the closure of $A$ and let $A_0$ be the
connected component of the identity of $\bar A$.  Let $V\subset N_{\mathbb{R}}$ be the
inverse image of $A_0$.  Then $A_0$ is a Lie group and so $V$ is a linear subspace.  As we
are assuming that $D$ is not rational, $A$ is infinite and so $A_0$ and $V$ are both
positive dimensional.  In particular $V$ is not contained in $\mathcal{G}$.  But then for
every $\epsilon>0$, we can find a positive multiple $mv$ of $v$, and a vector $w\in
N_{\mathbb{Z}}$, which is an integral linear combination of the divisors $P_k$, such that
\begin{itemize} 
\item $\ds \|mv-w\|<\epsilon$, whilst 
\item $mv-w\notin\mathcal{G}$.  
\end{itemize} 
Note that if $\epsilon>0$ is sufficiently small then $w\in \mathcal{P}$ since it is
integral and close to $mv\in\mathcal{P}$.  Thus if $M$ is the divisor corresponding to
$w$, then $M$ is base point free, and the rest is clear.  \end{proof}

\begin{definition}\label{d_saturated} Let $\pi\colon\map Y.Z.$ be a projective morphism
of normal varieties, where $Z$ is affine.  Let $\mathfrak{R}$ be the geometric algebra
associated to the additive sequence $M\bd$ of mobile divisors, with characteristic
sequence $D\bd$.
 
 We say that $\mathfrak{R}$ is \textbf{saturated} if there is a $\mathbb{Q}$-divisor 
$F$, such that 
\begin{enumerate} 
\item $\rup F.\geq 0$, and 
\item for every pair of positive integers $i$ and $j$,
$$
\mov (\rup jD_i+F.)\leq M_j.
$$
\end{enumerate} 
\end{definition}

\begin{theorem}\label{t_finite} Let $Y$ be a smooth variety and $\pi\colon\map Y.Z.$ a 
projective morphism, where $Z$ is affine and normal.  Let $\mathfrak{R}$ be a saturated
and free geometric ring on $Y$ whose characteristic sequence tends to a semiample limit.

 Then $\mathfrak{R}$ is finitely generated.  
\end{theorem}
\begin{proof} Let $D\bd$ be the characteristic sequence, with limit $D$.  Let 
$G$ be the support of $D$, and pick $\epsilon>0$ such that $\rup F-\epsilon G.\geq 0$.  
By diophantine approximation, we know that there is a positive integer $m$ and 
an integral divisor $M$ such that 
\begin{enumerate} 
\item $M$ is mobile, 
\item $\ds \|mD-M\| < \epsilon$, and 
\item If $mD\geq M$ then $mD=M$.  
\end{enumerate}

 But then
\begin{align*} 
mD+F &= M+(mD-M)+F\\ 
     &\geq M+F-\epsilon G, \\ 
\end{align*} 
so that 
$$
\mov (\rup mD+F.)\geq M.
$$ 
On the other hand, by definition of saturation we have
$$
\mov ( \rup mD_i+F.)\leq M_m=mD_m.
$$
Letting $i$ go to infinity we have
$$
M\leq \mov (\rup mD+F.)\leq mD_m\leq mD.
$$
By (3) above, it follows that the sequence of inequalities must in fact be equalities,
so that we have
$$
D=D_m,
$$
for some $m$ and we may apply \eqref{l_easy}.  
\end{proof}

\section{Saturation of the restricted algebra}
\label{s_saturation}

We fix some notation for this section.  Let $(X,\Delta)$ be a purely log terminal pair and
let $f\colon\map X.Z.$ be a projective morphism of normal varieties, where $Z$ is affine.
We assume that $S=\rdown\Delta.$ is irreducible.  Let $g\colon \map Y.X.$ be any log
resolution of the pair $(X,\Delta)$.  Then we may write
$$
K_Y+\Gamma=g^*(K_X+\Delta)+E,
$$
where $\Gamma$ and $E$ are effective, with no common components and $E$ is
$g$-exceptional.  We set $T=\rdown\Gamma.$ the strict transform of $S$, and
$F=E-\Gamma+T$.  We suppose that $K_Y+\Gamma$ is purely log terminal, so that $\rup F.\geq
0$ is effective and exceptional.  Fix a positive integer $k$ such that $k(K_X+\Delta)$ is
Cartier, so that both $G=k(K_Y+\Gamma)$ and $kE$ are integral.  Let $\pi\colon\map Y.Z.$
be the composition of $f$ and $g$.  Let $N_m+G_m$ be the decomposition of $mG$ into its
moving and fixed parts, and let $M_m$ be the restriction of $N_m$ to $T$.  Finally let
$$
D_i=\frac{M_i}i,
$$
so that $D\bd$ is the characteristic sequence of the restricted algebra
$\mathfrak{R}_T$.  

\begin{lemma}\label{l_ambient} For every pair of positive integers $i$ and $j$ 
$$
\mov (\rup (j/i)N_i+F.)\leq N_j.
$$
\end{lemma}
\begin{proof} We have, 
\begin{align*} 
\mov (\rup (j/i)N_i+F.) &\leq \mov (\rup (j/i)iG+F.) \\ 
                      &\leq \mov (\rup jk g^*(K_X+\Delta)+jkE+F.) \\ 
                      &\leq \mov (jk g^*(K_X+\Delta)+jkE+\rup F.) \\ 
                      &=\mov (jk g^*(K_X+\Delta)) \\ 
                      &\leq \mov (jk g^*(K_X+\Delta)+jkE) \\ 
                      &=\mov (jG) \\ 
                      &=N_j,\\ 
\end{align*} 
where we used the fact that $jkE+\rup F.$ is $g$-exceptional.  
\end{proof}

\begin{lemma}\label{l_saturation} Suppose that $-(K_X+\Delta)$ is nef and big
and that $M_m$ is free.

 Then the restricted algebra $\mathfrak{R}_T$ is saturated with respect to $F|_T$.  
\end{lemma}
\begin{proof} By assumption
$$
\rup F|_T.\geq 0.
$$

\begin{claim}\label{c_surjective} The natural restriction map,
$$
\map {H^0(Y,\ring Y.(\rup (j/i)N_i+F.))}.{H^0(T,\ring T.(\rup jD_i+F|_T.))}.,
$$
is surjective, for any positive integers $i$ and $j$.  
\end{claim}
\begin{proof}[Proof of claim] Fix $i$.  Since $M_i$ is free, it follows that $N_i$ is free
in a neighbourhood of $T$.  Then there is a model $h_i\colon\map Y_i.Y.$, on which $N_i$
becomes free, which is an isomorphism in a neighbourhood of $T$.  Thus, replacing $Y$ by
$Y_i$, we may assume that $N_i$ is free.
 
 Considering the restriction exact sequence,
$$
\ses {\ring Y.(\rup (j/i) N_i+F.-T)}.{\ring Y.(\rup (j/i)N_i+F.)}.{\ring T.(\rup jD_i+F|_T.)}.,
$$
it follows that the obstruction to the surjectivity of the restriction map above is given
by,
\begin{multline*}
H^1(Y,\ring Y.(\rup (j/i)N_i+(F-T).))\\
=H^1(Y,\ring Y.(K_{Y}+\rup g^*(-(K_X+\Delta)+(j/i)N_i).)),
\end{multline*}
which vanishes by Kawamata-Viehweg vanishing, as 
$$
g^*(-(K_X+\Delta)),
$$
is big and nef and $(j/i)N_i$ is nef.  
\end{proof}

 The result is now an easy consequence of \eqref{l_ambient} and the claim.  
\end{proof}

\begin{proof}[Proof of \eqref{t_restricted}] Let $g\colon\map Y.X.$ be the log resolution
of $(X,\Delta)$, whose existence is guaranteed by \eqref{t_reduce}.  We may write
$$
K_Y+\Gamma=g^*(K_X+\Delta)+E,
$$
where $\Gamma$ and $E$ are effective, with no common components, and $E$ is exceptional.
It follows that $k(K_Y+\Gamma)$ and $kE$ are integral and $mk(K_Y+\Gamma)$ and $mk
g^*(K_X+\Delta)$ have the same moving parts. Let $T$ be the strict transform of $S$ and
let $\mathfrak{R}_T$ be the restricted algebra associated to the divisor $k(K_Y+\Gamma)$.
By \eqref{l_birational}, it follows that $\mathfrak{R}_T$ is finitely generated iff
$\mathfrak{R}_S$ is finitely generated.  We may suppose that $\mathfrak{R}_T$ is a
limiting algebra, given by the convex sequence $\Theta\bd$.

By \eqref{c_universal} it follows that there is a birational model $h\colon\map T'.T.$,
where the mobile parts of $mk h^*(K_T+\Theta_m)$ are base point free, and that the limit
of the characteristic sequence is semiample.  We may assume that $h$ is induced by a
birational morphism $h'\colon\map Y'.Y.$.  Replacing $Y$ by $Y'$, we may assume that the
mobile parts of $mk(K_T+\Theta_m)$ are base point free and that the limit of the
characteristic sequence is semiample.

 It follows by \eqref{l_saturation} that the characteristic sequence is saturated.  
But then the restricted algebra $\mathfrak{R}_T$ is finitely generated by \eqref{t_finite}, and 
as we have already observed this implies that $\mathfrak{R}_S$ is finitely generated. 
\end{proof}

\begin{proof}[Proof of \eqref{t_existence}] Immediate from \eqref{t_restricted} and \eqref{l_implies}.  
\end{proof}

\begin{proof}[Proof of \eqref{c_existence1}] Clear.  \end{proof}

\begin{proof}[Proof of \eqref{c_existence2}] Follows from (5.1.3) of \cite{Shokurov96}.   \end{proof}

\bibliographystyle{hamsplain}
\bibliography{/home/mckernan/Jewel/Tex/math}

\end{document}